
\input amstex.tex
\documentstyle{amsppt}
\magnification1200
\hsize=12.5cm
\vsize=18cm
\hoffset=1cm
\voffset=2cm

\def\DJ{\leavevmode\setbox0=\hbox{D}\kern0pt\rlap
{\kern.04em\raise.188\ht0\hbox{-}}D}
\def\dj{\leavevmode
 \setbox0=\hbox{d}\kern0pt\rlap{\kern.215em\raise.46\ht0\hbox{-}}d}

\def\txt#1{{\textstyle{#1}}}
\baselineskip=13pt
\def\hf{{\textstyle{1\over2}}}
\def\a{\alpha}\def\b{\beta}
\def\d{{\,\roman d}}
\def\e{\varepsilon}
\def\f{\varphi}
\def\G{\Gamma}
\def\k{\kappa}
\def\s{\sigma}
\def\t{\theta}
\def\={\;=\;}
\def\zx{\zeta(\hf+ix)}
\def\zt{\zeta(\hf+it)}

\def\D{\Delta}
\def\no{\noindent}
\def\R{\Re{\roman e}\,} 
\def\z{\zeta}

 \def\t{\theta}
\def\hf{{\textstyle{1\over2}}}
\def\txt#1{{\textstyle{#1}}}
\def\f{\varphi}

\font\tenmsb=msbm10
\font\sevenmsb=msbm7
\font\fivemsb=msbm5
\newfam\msbfam
\textfont\msbfam=\tenmsb
\scriptfont\msbfam=\sevenmsb
\scriptscriptfont\msbfam=\fivemsb
\def\Bbb#1{{\fam\msbfam #1}}

\def \NN {\Bbb N}

\def \ZZ {\Bbb Z}

\font\ff=cmr8
\def\txt#1{{\textstyle{#1}}}
\baselineskip=13pt

\font\teneufm=eufm10
\font\seveneufm=eufm7
\font\fiveeufm=eufm5
\newfam\eufmfam
\textfont\eufmfam=\teneufm
\scriptfont\eufmfam=\seveneufm
\scriptscriptfont\eufmfam=\fiveeufm
\def\mathfrak#1{{\fam\eufmfam\relax#1}}

\font\tenmsb=msbm10
\font\sevenmsb=msbm7
\font\fivemsb=msbm5
\newfam\msbfam
     \textfont\msbfam=\tenmsb
      \scriptfont\msbfam=\sevenmsb
      \scriptscriptfont\msbfam=\fivemsb
\def\Bbb#1{{\fam\msbfam #1}}

\def \NN {\Bbb N}

\def \ZZ {\Bbb Z}

  \def\rightheadline{{\hfil{\ff
On the divisor function and the Riemann zeta-function in short intervals}
\hfil\tenrm\folio}}

  \def\leftheadline{{\tenrm\folio\hfil{\ff
   Aleksandar Ivi\'c }\hfil}}
  \def\emptyheadline{\hfil}
  \headline{\ifnum\pageno=1 \emptyheadline\else
  \ifodd\pageno \rightheadline \else \leftheadline\fi\fi}

\font\ff=cmr8
\font\teneufm=eufm10
\font\seveneufm=eufm7
\font\fiveeufm=eufm5
\newfam\eufmfam
\textfont\eufmfam=\teneufm
\scriptfont\eufmfam=\seveneufm
\scriptscriptfont\eufmfam=\fiveeufm
\def\mathfrak#1{{\fam\eufmfam\relax#1}}

\font\tenmsb=msbm10
\font\sevenmsb=msbm7
\font\fivemsb=msbm5
\newfam\msbfam
\textfont\msbfam=\tenmsb
\scriptfont\msbfam=\sevenmsb
\scriptscriptfont\msbfam=\fivemsb
\def\Bbb#1{{\fam\msbfam #1}}

\def \NN {\Bbb N}

\def \ZZ {\Bbb Z}

\def\D{\Delta}
\def\a{\alpha}
\def\b{\beta} \def\e{\varepsilon}
\def\no{\noindent} \def\d{\,{\roman d}}
\topmatter
\title
On the divisor function and the Riemann zeta-function in short intervals
\endtitle
\author
Aleksandar Ivi\'c
\endauthor
\address
Katedra Matematike RGF-a, Universitet u Beogradu,  \DJ u\v sina 7,
11000 Beograd, Serbia.
\endaddress
\subjclass 11 M 06, 11 N 37
\endsubjclass
\email
{\tt  ivic\@rgf.bg.ac.yu, aivic\@matf.bg.ac.yu}
\endemail
\dedicatory
\enddedicatory
\abstract
We obtain, for $T^\e \le U=U(T)\le T^{1/2-\e}$, asymptotic
formulas for
$$
\int_T^{2T}\left(E(t+U) - E(t)\right)^2\d t,\quad
\int_T^{2T}\left(\D(t+U) - \D(t)\right)^2\d t,
$$
where $\D(x)$ is the error term in the classical divisor problem,
and $E(T)$ is the error term in the mean square formula for $|\zt|$.
Upper bounds of the form $O_\e(T^{1+\e}U^2)$ for the above integrals
with biquadrates instead of square are shown to hold for
$T^{3/8} \le U =U(T) \ll T^{1/2}$. The connection between the moments of
$E(t+U) - E(t)$ and $|\zt|$ is also given. Generalizations to some other
number-theoretic error terms are discussed.
\endabstract
\endtopmatter
\document
\head
1. Introduction
\endhead
Power moments represent one of the most important parts of the
theory of the Riemann zeta-function $\z(s)$, defined as
$$
\z(s) = \sum_{n=1}^\infty n^{-s}\qquad(\s = \R s > 1),
$$
and otherwise by analytic continuation.
 Of particular significance are the moments on the
``critical line" $\s = \hf$, and a vast literature exists on this
subject (see e.g., the monographs [5], [6], and [23]). In this paper
we shall be concerned with moments of the error function
$$
E(T) := \int_0^T|\zt|^2\d t - T\left(\log{T\over2\pi} + 2\gamma -1\right),\leqno(1.1)
$$
where $\gamma = -\G'(1)$ is Euler's constant. More specifically, we shall consider
the moments
$$
\int_T^{2T}\left(E(t+G) - E(t-G)\right)^k\d t\qquad(k\in\NN\quad{\roman {fixed}}),
\leqno(1.2)
$$
where $G = G(T)$ is ``short" in the sense that $G =O(T)$ as
$T\to\infty$ and $G\gg1$. To deal with bounds for the expressions
like the one in (1.2), it seems convenient to use also results on
the moments of the function
$$
E^*(t) \;:=\; E(t) - 2\pi\D^*\bigl({t\over2\pi}\bigr),
$$
where
$$
\D^*(x) := -\D(x)  + 2\D(2x) - \hf\D(4x)
= \hf\sum_{n\le4x}(-1)^nd(n) - x(\log x + 2\gamma - 1).
$$
Here as usual $d(n) = \sum_{\delta|n}1$ is the number of positive
divisors of $n$, and
$$
\D(x) \;=\; \sum_{n\le x}d(n) - x(\log x + 2\gamma - 1)\leqno(1.3)
$$
is the error term in the classical Dirichlet divisor problem. The
function $E^*(t)$ gives an insight into the analogy between the
Dirichlet divisor problem and the mean square of $|\zt|$. It was
investigated by several authors, including M. Jutila [15], who
introduced the function $E^*(t)$, and the author [6]--[8].  Among
other things, the author (op. cit.) proved that
$$
\int_0^T (E^*(t))^2\d t \;=\; T^{4/3}P_3(\log T) + O_\e(T^{5/4+\e}),
$$
where $P_3$ is a polynomial of degree three in $\log T$ with
positive leading coefficient,
$$
\int_0^T |E^*(t)|^5\d t \;\ll_\e\; T^{2+\e},
\quad
\int_0^T|E^*(t)|^3 \d t \;\ll_\e\; T^{3/2+\e},\leqno(1.4)
$$
and none of these three results implies any one of the other two.
From the bounds in (1.4) and the Cauchy-Schwarz inequality for integrals it
follows that
$$
\int_0^T |E^*(t)|^4\d t \;\ll_\e\; T^{7/4+\e}.\leqno(1.5)
$$
Here and later $\e \;(>0)$ denotes arbitrarily small constants,
not necessarily the same ones at each occurrence, and $a = O_\e(b)$
(same as $a \ll_\e b$) means that the implied constant depends only on $\e$.
In addition to (1.2) it makes sense to investigate the moments
$$
\int_T^{2T}\left(\D(t+G) - \D(t-G)\right)^k\d t\qquad(k\in\NN\quad{\roman {fixed}}),
\leqno(1.6)
$$
as well. The interest in this topic comes from the work of M. Jutila
[12], who investigated the case $k=2$ in (1.2) and (1.6). He proved
that
$$
\eqalign{& \int\limits_T^{T+H}\left(\D(x+U)-\D(x)\right)^2\d x \cr&=
{1\over4\pi^2}\sum_{n\le {T\over2U}}{d^2(n)\over
n^{3/2}}\int\limits_T^{T+H} x^{1/2}\left|\exp\left(2\pi iU\sqrt{{n\over
x}}\,\right)-1\right|^2\d x + O_\e(T^{1+\e} +
HU^{1/2}T^\e),\cr}\leqno(1.7)
 $$
 for $1 \le U \ll T^{1/2} \ll H \le T$,
and an analogous result holds also for the integral of $E(x+U)-E(x)$
(the constants in front of the sum and in the exponential will be
$1/\sqrt{2\pi}$ and $\sqrt{2\pi}$, respectively).
From (1.7) one deduces ($a\asymp b$ means  $a\ll b\ll a$)
$$
\int_T^{T+H}\left(\D(x+U)-\D(x)\right)^2\d x
\asymp HU\log^3\left({\sqrt{T}\over U}\right)
\leqno(1.8)
$$
for $HU \gg T^{1+\e}$ and $T^\e \ll U \le \hf\sqrt{T}$. In [14]
Jutila proved that the integral in (1.8) is
$$
\ll_\e T^\e(HU + T^{2/3}U^{4/3})\qquad(1 \ll H,U \ll X).
$$
This bound and (1.8) hold also for the integral of $E(x+U)-E(x)$.
Furthermore Jutila conjectured that
$$
\int_T^{2T}\left(E(t+U) - E(t-U)\right)^4\d t \ll_\e T^{1+\e}U^2\leqno(1.9)
$$
holds for $1 \ll U \ll T^{1/2}$, and the analogous formula should hold for
$\D(t)$ as well. In fact, using the ideas of K.-M. Tsang [24] who
investigated the fourth moment of $\D(x)$, it can be shown that one
expects the integral in (1.9) to be of order
$TU^2\log^6(\sqrt{T}/U)$. Jutila also indicated that the truth of
his conjecture (1.9) implies
$$
\int_0^T|\zt|^6\d t \;\ll_\e\;T^{1+\e}.\leqno(1.10)
$$
This is (a weakened form of) the sixth moment for $|\zt|$, and the
best known exponent at present on the right-hand side of (1.10) is
5/4 (see [5], [6]). In view of the bound (op. cit.)
$$
|\zt|^k \ll \log t\int_{t-1}^{t+1}|\zx|^k \d x + 1,
\qquad(k\in \NN\;\;{\roman{fixed}})
\leqno(1.11)
$$
we actually have, using (1.9) with $U = T^\e$ and (1.11) with $k=2$,
$$
\int_T^{2T}|\zt|^8\d t \ll_\e \int_T^{2T}\left\{\log T(E(t+T^\e)-E(t-T^\e))^4+T^\e
\right\}\d t \ll_\e T^{1+\e},\leqno(1.12)
$$
and the eighth moment bound (1.12) is notably stronger than (1.10).
It may be remarked that the fourth moments of $\D(x)$ and $E(T)$
have been investigated by several authors, including Ivi\'c--Sargos
[11], K.-M. Tsang [24], and W. Zhai [25], [26].

\bigskip
\head
2. Statement of results
\endhead
\no Our first aim is to derive from (1.7) (when $H=T$)  a true asymptotic
formula. The result is

\bigskip
THEOREM 1. {\it For $1 \ll U = U(T) \le \hf {\sqrt{T}}$ we have} ($c_3 = 8\pi^{-2}$)
$$\eqalign{
\int_T^{2T}\left(\D(x+U)-\D(x)\right)^2\d x & = TU\sum_{j=0}^3c_j\log^j
\Bigl({\sqrt{T}\over U}\Bigr) \cr&
+ O_\e(T^{1/2+\e}U^2) + O_\e(T^{1+\e}U^{1/2}),\cr}\leqno(2.1)
$$
{\it a similar result being true
if $\D(x+U)-\D(x)$ is replaced by $E(x+U)-E(x)$, with different constants $c_j$.}

\medskip
{\bf Remark 1}. For $T^\e \le U = U(T) \le T^{1/2-\e}$ (2.1) is a true asymptotic formula.

\medskip
{\bf Corollary 1}. For $1 \ll U \le \hf {\sqrt{T}}$ we have ($c_3 = 8\pi^{-2}$)
$$\eqalign{
\sum_{T\le n \le 2T}\left(\D(n+U)-\D(n)\right)^2 & = TU\sum_{j=0}^3c_j\log^j
\Bigl({\sqrt{T}\over U}\Bigr) \cr&
+ O_\e(T^{1/2+\e}U^2) + O_\e(T^{1+\e}U^{1/2}),\cr}\leqno(2.2)
$$
The formula (2.2) is a considerable improvement over a result of Coppola--Salerno [3],
who had ($T^\e \le U \le \hf \sqrt{T},\,L = \log T$)
$$
\sum_{T\le n \le 2T}\left(\D(n+U)-\D(n)\right)^2=
{8\over\pi^2}TU\log^3\Bigl({\sqrt{T}\over U}\Bigr)
+ O(TUL^{5/2}\sqrt{L}).\leqno(2.3)
$$

\medskip
{\bf Corollary 2}. For $T \le x \le 2T$ and $T^\e \le U= U(T) \le T^{1/2-\e}$ we have
$$
\D(x+U) - \D(x) = \Omega \Bigl\{\sqrt{U}\log^{3/2}\Bigl({\sqrt{x}\over U}\Bigr)\Bigr\},\;
E(x+h) - E(x) = \Omega \Bigl\{\sqrt{U}\log^{3/2}\Bigl({\sqrt{x}\over U}\Bigr)\Bigr\}.
\leqno(2.4)
$$
These omega results ($f(x) = \Omega(g(x))$
means that $\lim_{x\to\infty}f(x)/g(x)\ne0$) show that Jutila's conjectures made in [12],
namely that
$$
\D(x+U) - \D(x) \ll_\e x^\e\sqrt{U},\;E(x+U) - E(x) \ll_\e x^\e\sqrt{U}\leqno(2.5)
$$
for $x^\e \le U \le x^{1/2-\e}$ are (if true), close to being best possible.
The difficulty of these conjectures may be seen if one notes that from the
definition of $\D(x)$ (the analogue of this for $E(T)$ is not known to hold,
in fact it is equivalent to the Lindel\"of hypothesis (see [6])) one easily obtains
$$
\D(x+U) - \D(x) \ll_\e x^\e{U}\qquad(1 \ll U \le x),\leqno(2.6)
$$
which is much weaker than (2.5). However,
a proof of (2.6) has not been obtained yet by the classical Vorono{\"\i} formula.
This formula will be needed later for the proof of Theorem 3,
and in a  truncated form it reads (see e.g., Chapter 3 of [5])
$$
\D(x) = {1\over\pi\sqrt{2}}
x^{1\over4}\sum_{n\le N}d(n)n^{-{3\over4}}\cos(4\pi\sqrt{nx}
- {\txt{1\over4}}\pi) +
O_\e(x^{{1\over2}+\e}N^{-{1\over2}})\quad(2 \le N \ll x).
\leqno(2.7)
$$
One also has (see  [5, eq. (15.68)]), for $2 \le N \ll x$,
$$
\D^*(x) = {1\over\pi\sqrt{2}}x^{1\over4}
\sum_{n\le N}(-1)^nd(n)n^{-{3\over4}}
\cos(4\pi\sqrt{nx} - {\txt{1\over4}}\pi) +
O_\e(x^{{1\over2}+\e}N^{-{1\over2}}),
\leqno(2.8)
$$
which is completely analogous to (2.7), the only difference is that in (2.8)
there appears a factor $(-1)^n$ in the sum.

\medskip
{\bf Remark 2}. The analogue of (2.3) for the sum
$$
\sum_{T\le n \le 2T}\left(E(n+U)-E(n)\right)^2\leqno(2.9)
$$
does not carry over, because $E(T)$ (see (1.1)) is a continuous function,
while $\D(x)$ is not, having jumps at natural numbers of order at most
$O_\e(x^\e)$. The true order of magnitude of the sum in (2.9) seems elusive.
Sums of $E(n)$ were investigated by Y. Bugeaud and the author [2]. By
using the irrationality measure of
${\roman e}^{2\pi m}$ and  for the partial quotients
in its continued fraction expansion, a non-trivial bound
for $\sum_{n\le x}E(n)$ is obtained.

\medskip
There are several ways in which the asymptotic formula (2.1) of Theorem 1
may be generalized. This concerns primarily number-theoretic terms related
to arithmetic functions $f(n)$ whose generating series $F(s) =
\sum_{n=1}^\infty f(n)n^{-s}\;(\Re s >1)$ belongs to the so-called
Selberg class of degree two (see e.g., the survey work of Kaczorowski--Perelli
[16]). Instead of trying to formulate a general result which contains (2.1)
as a special case, we shall state the corresponding results for two well-known
number-theoretic quantities. Let, as usual, $r(n) = \sum_{n=a^2+b^2}1$ denote
the number of ways $n$ may be represented as a sum of two integer squares, and
let  $\varphi(z)$ be a
holomorphic cusp form of weight $\kappa$ with respect to the full
modular group $SL(2,\ZZ)$, and denote by $a(n)$ the $n$-th Fourier
coefficient of $\varphi(z)$. We suppose that $\varphi(z)$ is a
normalized eigenfunction for the Hecke operators $T(n)$, that is,
$  a(1)=1  $ and $  T(n)\varphi=a(n)\varphi $ for every $n \in
\NN$ (see e.g., R.A. Rankin [20] for the definition and properties
of the Hecke operators). The classical example is $a(n) = \tau(n)\;(\k=12)$,
 the Ramanujan function defined by
$$
\sum_{n=1}^\infty \tau(n)x^n \=
x{\left\{(1-x)(1-x^2)(1-x^3)\cdots\right\}}^{24}\qquad(\,|x| < 1).
$$
If $P(x) := \sum_{n\le x }r(n) - \pi x$ denotes then the error term in
the classical circle problem and $A(x) := \sum_{n\le x}a(n)$, then we have

\medskip
THEOREM 2. {\it For $T^\e \le U=U(T) \le  \hf\sqrt{T}$ we have}
$$\eqalign{
\int_T^{2T}(P(t+U)-P(t))^2\d t &= TU\left(A_1\log\Bigl({\sqrt{T}\over U}\Bigr)
+ A_2\right)\cr&
+ O_\e(T^{1/2+\e}U^2) + O_\e(T^{1+\e}\sqrt{U}),\cr}
\leqno(2.10)
$$
{\it and}
$$
\int_T^{2T}(A(t+U)-A(t))^2\d t = CT^\kappa U
+ O_\e(T^{\kappa-2/5+\e}U^{9/5}) + O_\e(T^{\kappa+\e}\sqrt{U})
\leqno(2.11)
$$
{\it with some explicitly computable constants $A_1,C >0$ and $A_2$.}

\medskip
{\bf Corollary 3}. For $T \le x \le 2T$ and $T^\e \le U= U(T) \le T^{1/2-\e}$ we have
$$
P(x+U) - P(x) = \Omega\left(\sqrt{U\log\Bigl({\sqrt{x}\over U}\Bigr)}\,\right),
\quad A(x+U)-A(x) = \Omega(x^{{\kappa-1\over2}}\sqrt{U}).
$$
\medskip
Our next  result relates bounds for
moments of $|\zt|$ to bounds of moments of $E(t+G) - E(t-G)$. This
is usually done (see e.g., Chapter 8 of [5]) by counting ``large
values" of $|\zt|$ which occur in $[T, 2T]$. Our result is the
following
\bigskip
THEOREM 3. {\it Let $t_1,\ldots, t_R$ be points in $[T,\,2T]$ which
satisfy $T^\e \le V \le  |\z(\hf + it_r)|$  and $|t_r - t_s|\ge 1$ for $r,s \le R$
and $r\ne s$. Then we have, for $L = \log T, G = A(V/L)^2$ with a
 suitable constant $A>0$,and $k\in \NN$ fixed,}
$$
R \ll V^{-2-2k}L^{2+2k}\int\limits_{T/3}^{3T}\Bigl\{|E(t+2G)-E(t-2G)|^k
+|E(t+\hf G)-E(t-\hf G)|^k\Bigr\}\d t.\leqno(2.12)
$$

\bigskip
{\bf Corollary 4}. {\it Suppose that the integral on the right-hand side of} (2.12)
{\it is bounded by $O_\e(T^{\a+\e}G^\b)$ for some real constants
$\a = \a(k)\;(>0)$ and $\b = \b(k) \le k-1$, and $T^\e \le G = G(T) \ll T^{1/3}$.
Then we have}
$$
\int_0^T|\zt|^{2+2k-2\b}\d t \;\ll_\e\;T^{1+\a+\e}.\leqno(2.13)
$$

\bigskip
One obtains  Corollary 3 from Theorem 3 in a standard  way (see
e.g., Chapter 8 of [5]). The condition $T^\e \le G \ll T^{1/3}$
comes from the definition of $G$ and the classical bound $\zt \ll
t^{1/6}$. The condition $\b \le k-1$ is necessary, because we know
that $\int_0^T|\zt|^4\d t \ll T\log^4T$, and the condition in
question implies that the exponent of the integral in (2.13) is at
least 4.

\medskip
In connection with Jutila's conjecture (1.9) one may, in general, consider
constants $0 \le \rho(k) \le 1$ for fixed $k>2$ which one has
$$
\int_T^{2T}|E(t+G)-E(t-G)|^k\d t \ll_\e
T^{1+\e}G^{k/2}\quad(T^{\rho(k)+\e} \ll G = G(T) \ll T),
\leqno(2.14)
$$
and similarly for the moments of $|\D(t+G)-\D(t-G)|$. A general, sharp version of Jutila's
conjecture would be that $\rho(k)=0$ for any fixed $k>2$ and (2.14) holds for $G\ll \sqrt{T}$.
The following theorem gives the unconditional value of $\rho(4)$, and
shows that Jutila's conjecture holds in a certain range. Any improvements of these
ranges would be of interest.

\bigskip
THEOREM 4. {\it We have, for  $\;T^{3/8} \ll G = G(T) \ll T^{1/2}$,}
$$\eqalign{&
\int\limits_T^{2T}(E(t+G)-E(t-G))^4\d t \;\ll_\e\; T^{1+\e}G^2,\cr&
\int\limits_T^{2T}(\D(t+G)-\D(t-G))^4\d t \;\ll_\e\; T^{1+\e}G^2.\cr}
\leqno(2.15)
$$

\break
\head
3. The proof of Theorem 1 and Theorem 2
\endhead
We shall deduce Theorem 1 from Jutila's formula
(1.7) with $H = T$. First note that
the integral on the right-hand side equals
$$
\eqalign{&
\int_T^{2T}x^{1/2}\left|\exp\left(2\pi iU\sqrt{{n\over x}}\,\right)
-1\right|^2\d x\cr&
= \int_T^{2T}x^{1/2}\Bigl(
2- {\roman e}^{-2\pi iU\sqrt{n/x}} - {\roman e}^{2\pi iU\sqrt{n/x}}\,\Bigr)\d x\cr&
= 2\int_T^{2T}x^{1/2}\left(1-\cos\Biggl(2\pi U\sqrt{{n\over x}}\,\Biggr)\right)\d x
\cr&
= 4\int_T^{2T}x^{1/2}\sin^2\Biggl(\pi U\sqrt{{n\over x}}\,\Biggr)\d x.\cr}
$$
In the last integral we make the change of variable
$$
\pi U\sqrt{{n\over x}} = y,\; \sqrt{x} = {\pi U\sqrt{n}\over y},\; x = \pi^2U^2ny^{-2},
\;\d x = -2\pi^2U^2ny^{-3}.
$$
Therefore the main term on the right-hand side of (1.7) becomes
$$
2\pi U^3\sum_{n\le T/(2U)}d^2(n)\int_{\pi U\sqrt{n/(2T)}}^{\pi U\sqrt{n/T}}
{\sin^2y\over y^4}\,\d y.\leqno(3.1)
$$
Now we change the order of summation and integration: from
$$
1 \le n\le {T\over2U},\; \pi U\sqrt{{n\over2T}} \le y \le \pi U\sqrt{{n\over T}}
$$
we infer that
$$
{\pi U\over\sqrt{2T}} \le y \le \pi\sqrt{U\over2},\quad {Ty^2\over \pi^2 U^2}
\le n \le {2Ty^2\over \pi^2 U^2}.
$$
Thus (3.1) becomes
$$
2\pi U^3\int_{\pi U\over\sqrt{2T}}^{\pi\sqrt{U\over2}}
\sum_{\max(1,{Ty^2\over \pi U^2})\le n \le \min({T\over2U},{2Ty^2\over \pi U^2})}d^2(n)
\cdot {\sin^2y\over y^4}\,\d y.\leqno(3.2)
$$
The range of summation in (3.2) will be
$$
I\;:=\;\left[{Ty^2\over \pi U^2}, {2Ty^2\over \pi U^2}\right],\quad
{\roman {if}}\;\, y\in J:= \left[{\pi U\over\sqrt{T}},\hf\pi\sqrt{U}\,\right].
$$
By using the elementary bound
$|\sin x| \le \min(1,|x|)$, it is easily seen that the error made by replacing the interval
of integration in (3.2) by $J$ will be
$$
\ll (TU + T^{1/2}U^2)\log^3T,
$$
which is absorbed by the error term in (2.1).
When $y\in J$, the sum over $n\in I$ can be evaluated by the use of the asymptotic formula
(see [5] and [10, Lemma 3])
$$
\sum_{n\le x}d^2(n) = x\Bigl(\sum_{j=0}^3a_j\log^jx\Bigr) + O_\e(x^{1/2+\e})
\quad(a_3 = 1/(\pi^2)).\leqno(3.3)
$$
We note that the value $a_3 = 1/(\pi^2)$ is easily computed from the residue
of $x^s\z^4(s)/s\z(2s)$ at $s=1$, and the remaining $a_j$'s in (3.3) can be also
explicitly computed. The error term in (3.3) can be improved
to $O(x^{1/2}\log^5x\log\log x)$ (see Ramachandra--Sankaranarayanan [18]),
but the exponent $1/2$ of $x$ cannot be improved without
assumptions on the zero-free region of $\z(s)$ (such as e.g., the Riemann
hypothesis that all complex zeros of $\z(s)$ have real parts equal to 1/2).
However, this improvement is not needed in view of the error term $O(HU^{1/2}T^\e)$
in (1.7).

\smallskip
To continue with the proof, note that if we use (3.3) to
evaluate the expression in (3.2) we shall obtain,
with effectively computable constants $b_j\, (b_3 = 1/(\pi^2))$, that
the major contribution equals
$$
\eqalign{&
2\pi U^3\int_{\pi U\over\sqrt{T}}^{{1\over2}\pi\sqrt{U}}{\sin^2y\over y^4}
\Biggl\{{Ty^2\over\pi^2U^2}\Biggl(\sum_{j=0}^3b_j\log^j\Bigl({Ty^2\over U^2}\Bigr)\Biggr)
+ O_\e\Bigl({T^{1/2+\e}y\over U}\Bigr)\Biggr\}\d y\cr&
= {2\over\pi}TU\int_{\pi U\over\sqrt{T}}^{{1\over2}\pi\sqrt{U}}{\sin^2y\over y^2}
\Biggl(\sum_{j=0}^3b_j\log^j\Bigl({Ty^2\over U^2}\Bigr)\Biggr)\d y + O_\e(T^{1/2+\e}U^{2}).
\cr}\leqno(3.4)
$$
The last error term above comes from the fact that
$$
\eqalign{
\int_{\pi U\over\sqrt{T}}^{{1\over2}\pi\sqrt{U}}{\sin^2y\over y^3}\d y&
= \int_{\pi U\over\sqrt{T}}^{1}{\sin^2y\over y^3}\d y + O(1)\cr&
\ll \int_{\pi U\over\sqrt{T}}^{1}{\d y\over y} + 1 \ll \log {\sqrt{T}\over U},
\cr}
$$
where  $|\sin x| \le \min(1,|x|)$ was used again. Likewise we deduce that,
for $0<\a \le 1, \,\b\gg1$,
$$
\eqalign{
\int_\a^\b\,{\sin^2y\over y^2}\d y &= \int_0^\infty\,{\sin^2y\over y^2}\d y + O(\a) + O(\b^{-1})
\cr&
= {\pi\over2} + O(\a) + O(\b^{-1}).\cr}\leqno(3.5)
$$
We expand as a binomial
$$
\log^j\Bigl({Ty^2\over U^2}\Bigr) = {\Bigl(\log{T\over U^2}+ 2\log y\Bigr)}^j
\qquad(j = 2,3),
$$
and use a relation similar to (3.5) for an integral
containing an additional power of $\log y$.
Hence from (3.4) it transpires that the main term on the right-hand side of (1.7) is
equal to
$$
\eqalign{&
{2\over\pi}TU\Biggl\{{\pi\over2}b_3\log^3\bigl({T\over U^2}\bigr) +
c_2'\log^2\bigl({T\over U^2}\bigr) + c_1'\log\bigl({T\over U^2}\bigr) + c_0'\cr&
+ O_\e(T^{\e-1/2}U + T^\e U^{-1/2})\Biggr\}\cr&
= TU\left\{{8\over\pi^2}\log^3\Bigl({\sqrt{T}\over U}\Bigr)
+ c_2\log^2\Bigl({\sqrt{T}\over U}\Bigr) + c_1\log\Bigl({\sqrt{T}\over U}\Bigr)
+ c_0\right\}\cr&
 + O_\e(T^{1/2+\e}U^2 + T^{1+\e}U^{1/2}).
\cr}\leqno(3.6)
$$
From (3.6) and (1.7) we easily obtain (2.1). The proof of (2.1) with $E(x+U)-E(x)$
in place of $\D(x+U)-\D(x)$ follows verbatim the above argument.

\medskip
The formula (2.1) of Corollary 1 follows from (2.1) and
$$
\int_T^{2T}(\D(x+U)-\D(x))^2\d x = \sum_{T\le n\le 2T}(\D(n+U)-\D(n))^2
+ O(U^{5/2}\log^{5/2}T),\leqno(3.7)
$$
for $1 \ll U \ll \sqrt{T}$. Namely we can assume $U,T$ are integers (otherwise
making an admissible error). Using (1.3) and the mean
value theorem, it follows that the left-hand side
of (3.7) equals ($0\le \t\le1$)
$$
\eqalign{&
\sum_{T\le m\le 2T-1}\int_m^{m+1-0}\left(\sum_{x<n\le x+U}d(n) - U\Bigl
(\log(x+\t U)+2\gamma\Bigr)\right)^2\d x
\cr&
= \sum_{T\le m\le 2T-1}\int_m^{m+1-0}\left(\sum_{m<n\le m+U}d(n) -
U\Bigl(\log(x+\t U)+2\gamma\Bigr)\right)^2\d x
\cr&
= \sum_{T\le m\le 2T-1}\left(\D(m+U) - \D(m) + O(U^2T^{-1}\log T)\right)^2.\cr}
$$
Now we expand the square, use the Cauchy-Schwarz inequality and (2.3)
for the cross terms, replace the range of summation by $[T,\,2T]\,$, and (3.7) follows.

\medskip
To prove Theorem 2, note first that $r(n) = 4\sum_{d|n}\chi(d)$, where $\chi$ is the
non-principal character modulo four. Thus $0 \le r(n) \le 4d(n)$, and ${1\over4}r(n)$
is multiplicative. From the functional equation
$$
L(s) = \pi^{2s-1}{\G(1-s)\over\G(s)}L(1-s),
$$
where $L(s)$ is the generating Dirichlet series of $r(n)$ one obtains the explicit formula
$$
P(x) = -{1\over\pi}x^{1/4}\sum_{n\le N}r(n)n^{-3/4}\cos(2\pi\sqrt{nx} + {\pi\over4})
+ O_\e(x^{1/2+\e}N^{-1/2})\leqno(3.8)
$$
for $1\ll N\ll x$, much in the same way as one obtains (2.7) (see e.g., Chapter
13 of [5]). This formula is completely analogous to (2.7), and consequently
Jutila's proof gives the analogue of (1.7), with a different constant in front of the
sum, $d(n)$ replaced by $r(n)$, and $\pi i U\sqrt{n/x}$ in the exponential.
The proof of Theorem 1 goes through up to (3.3), where instead of the asymptotic
formula for sums of $d^2(n)$ we shall use
$$
\sum_{n\le x}r^2(n) = 4x\log x + Cx + O(x^{1/2}\log^3x\log\log x),\leqno(3.9)
$$
where $C = 8.0665\ldots$ is an explicitly given constant.
The asymptotic formula (3.9) is due to K\"uhleitner--Nowak
[17]. Note that the main term in (3.9) is somewhat different than the main term in (3.3),
which is reflected in different main terms in (2.1) and (2.10). By using (3.9) the proof
of (2.10) is essentially the same as the proof of (2.1), so there is no need
for the details. For our purposes (3.9) with the error term $O_\e(x^{12+\e})$
suffices, in view of the term $O_\e(T^{1+\e}U^{1/2})$ in the analogue of (2.1).

\medskip
As to the proof of (2.11), note that by P. Deligne's bound one has
$|a(n)| \le n^{(\kappa-1)/2}d(n)$, and the analogue of (2.7) reads
(see e.g., M. Jutila [13] for a proof, who has a more general result with
exponential factors)
$$\eqalign{
A(x) = \sum_{n\le x}a(n) &={1\over\pi\sqrt{2}}x^{{\kappa\over4}-{1\over4}}
\sum_{n\le N}a(n)
n^{-{\kappa\over2}-{1\over4}}\cos\Bigl(4\pi\sqrt{nx}-{\pi\over4}\Bigr)\cr&
+ O_\e(x^{\kappa/2+\e}n^{-1/2})\qquad(1\ll N \ll x).\cr}\leqno(3.10)
$$
One has the asymptotic formula ($A>0$ can be explicitly evaluated)
$$
\sum_{n\le x}a^2(n) = Ax^\kappa + O(x^{\kappa-2/5}).\leqno(3.11)
$$
The above formulas show that $a(n)$ behaves similarly to $n^{(\kappa-1)/2}d(n)$.
The bound for the error term in (3.11), one of the longest standing records
in analytic number theory is due to R.A. Rankin  [19] and A. Selberg [22].
Following Jutila's proof of (2.1), the proof of (2.10) and using (3.10) instead
of (3.3) at the appropriate place, we arrive at (2.11).

\head
4. The proof of Theorem 3
\endhead
In this section we shall present the proof of Theorem 3.
From the definition (1.1) of $E(T)$ we have, for $T\le u,t\le 2T,\,1\ll G\ll T$,
$$
E(u+\hf G) - E(u-\hf G) = \int_{u-G/2}^{u+G/2}|\zx|^2\d x + O(G\log T).
$$
Consequently integration over $u$ gives
$$
\eqalign{&
\int_{t-G/2}^{t+G/2}(E(u+\hf G) - E(u-\hf G))\d u
\cr& = \int_{t-G/2}^{t+G/2}\int_{u-G/2}^{u+G/2}|\zx|^2\d x \d u + O(G^2\log T)\cr&
\le \int_{t-G/2}^{t+G/2}\int_{t-G}^{t+G}|\zx|^2\d x\d u + O(G^2\log T)\cr&
= G\int_{t-G}^{t+G}|\zx|^2\d x  + O(G^2\log T).\cr}
$$
Using again (1.1) for the last integral it follows that
$$
E(t+G) - E(t-G) \ge {1\over G}\int_{t-G/2}^{t+G/2}(E(u+\hf G) - E(u-\hf G))\d u
- CG\log T\leqno(4.1)
$$
for $1\ll G \ll T$ and a suitable constant $C>0$.
The bound in (4.1) is useful when $E(t+G) - E(t-G)$ is negative.
Likewise, from
$$
\eqalign{&
\int_{t-G}^{t+G}(E(u+ 2G) - E(u-2G))\d u \cr&
= \int_{t-G}^{t+G}\int_{u-2G}^{u+2G}|\zx|^2\d x\d u + O(G^2\log T)\cr&
\ge \int_{t-G}^{t+G}\int_{t-G}^{t+G}|\zx|^2\d x\d u + O(G^2\log T)\cr&
= 2G\int_{t-G}^{t+G}|\zx|^2\d x + O(G^2\log T)\cr&
= 2G(E(t+G) - E(t-G)) + O(G^2\log T)\cr}
$$
we obtain a bound which is useful when $E(t+G) - E(t-G)$ is positive. This is
$$
E(t+G) - E(t-G) \le {1\over 2G}\int_{t-G}^{t+G}(E(u+2G) - E(u-2G))\d u
+ CG\log T.\leqno(4.2)
$$
Combining (4.1) and (4.2), depending on the sign of $E(t+G) - E(t-G)$,
we obtain, for $T \le t \le 2T,\,1\ll G\ll T,\,C>0$,
$$
\eqalign{&
|E(t+G) - E(t-G)| \le CG\log T +\cr&
+ {1\over G}\int_{t-G}^{t+G}\Bigl\{|E(u+2G)-E(u-2G)| + |E(u+\hf G)-E(u-\hf G)|\Bigr\}\d u.
\cr}\leqno(4.3)
$$
Suppose now that the hypotheses of Theorem 3 hold. Then ($L = \log T)$
$$
V^2 \le |\z(\hf+it_r)|^2 \ll L\left(\int_{t_r-1/3}^{t_r+1/3}|\zx|^2\d x+1\right)
\qquad(r =1,\ldots\,,R).\leqno(4.4)
$$
The interval $[T,\,2T]$ is covered then with subintervals of length $2G$, of which
the last one may be shorter. In these intervals we group subintegrals over disjoint
intervals $[t_r-1/3,\,t_r+1/3]$. Should some intervals fall into two of such
intervals of length $2G$, they are treated then separately in an analogous
manner. It follows that
$$
R \ll V^{-2}L^2\sum_{j=1}^J \int_{\tau_j-G}^{\tau_j-G}|\zx|^2\d x,
$$
where $J\le R$, $\tau_j \in [T/3,3T],\, |\tau_j-\tau_\ell|\ge 2G\;
(j\ne \ell; j,\ell \le J)$
by considering separately points with even and odd indices.
Now we note that by (1.1)
$$
\int_{\tau_j-G}^{\tau_j-G}|\zx|^2\d x = O(GL) + E(\tau_j-G)-E(\tau_j+G).\leqno(4.5)
$$
For $E(\tau_j-G)-E(\tau_j+G)$ we use (4.3) with $t = \tau_j$, choosing
$$
G = AV^2L^{-2}
$$
with suitable $A (> 0)$ so that $O(GL) \le \hf V^2$.  In this way we obtain,
using H\"older's inequality for integrals, noting that the intervals
$[\tau_j- G,\,\tau_j +G]$ are disjoint (if we consider separately systems
of points $\tau_j$ with even and odd indices $j$) and $J\le R$,
$$
\eqalign{&
R \ll V^{-4}L^4\sum_{j=1}^J \int_{\tau_j-G}^{\tau_j-G}\left\{|E(u+2G) -\cdots|\right\}\d u
\cr&
\ll V^{-4}L^4\sum_{j=1}^J\left(\int_{\tau_j-G}^{\tau_j-G}\left\{|E(u+2G) - \cdots|
\right\}^k\d u\right)^{1/k}G^{1-1/k}\cr&
\ll V^{-4}L^4(RG)^{1-1/k}\left(\int_{T/3}^{3T}
\left\{|E(u+2G)- \cdots|\right\}^k\d u\right)^{1/k}.
\cr} \leqno(4.6)
$$
If we simplify (4.6), we obtain the assertion (2.12) of Theorem 3.
\medskip
\break

\head
5. The proof of Theorem 4
\endhead

For the proof of Theorem 4 we shall need the case $k=2$ of the following
\medskip
LEMMA 1.  {\it Let $k\ge 2$ be a fixed
integer and $\delta > 0$ be given.
Then the number of integers $n_1,n_2,n_3,n_4$ such that
$N < n_1,n_2,n_3,n_4 \le 2N$ and}
$$
|n_1^{1/k} + n_2^{1/k} - n_3^{1/k} - n_4^{1/k}| < \delta N^{1/k}
$$
{\it is, for any given $\e>0$,}
$$
\ll_\e N^\e(N^4\delta + N^2).\leqno(5.1)
$$
\medskip
Lemma 1 was proved by Robert--Sargos [21]. It represents a powerful
arithmetic tool which is essential in the analysis when the
biquadrate of sums involving $\sqrt n$ appears in exponentials, and
was used e.g., in [11].

\medskip
It is enough to prove (2.14) of Theorem 4 for $\D(x)$. Namely because of the
analogy between (2.7) and (2.8) (which differs from (2.13) only by the presence
of the innocuous factor $(-1)^n$ in the sum), the same bound in the same range
for $G$ will hold with the integral of $\D^*(x)$ replacing $\D(x)$. But then,
in view of
$$
E(t) = E^*(t) +  2\pi\D^*\bigl({t\over2\pi}\bigr)
$$
and (1.5), we obtain
$$\eqalign{&
\int_T^{2T}(E(t+G) - E(t-G))^4\d t \cr&
\ll \int_T^{2T}\Bigl(E^*(t+G) - E^*(t-G)\Bigr)^4\d t
+ \int_T^{2T}\left(\D^*({t+G\over2\pi}) - \D^*({t-G\over2\pi})\right)^4\d t\cr&
\ll_\e T^{7/4+\e} + T^{1+\e}G^2\ll_\e T^{1+\e}G^2\cr}
$$
precisely for $G\ge T^{3/8}$.

\bigskip
For the proof of (2.14) with $\D(x)$ we start from (2.7) with $x = t+G$,
$x = t-G, T\le t \le 2T,
N = T$ in both cases. We split the sum over $n$ into $O(\log T)$
subsums over $M < n\le M'\le 2M$,
and raise each sum in question to the fourth power and integrate.
When $M \ge TG^{-4/3}$ we note that, using twice (2.7), we have
$$
S(t,M) := t^{1/4}\sum_{M<n\le M'}d(n)n^{-3/4}\cos(4\pi\sqrt{nt}-\pi/4)
\ll_\e T^{1/2+\e}M^{-1/2}. \leqno(5.2)
$$
Hence from (5.2) and the first derivative test (see e.g., Lemma 2.1
of [5]) we infer that in this range
$$
\eqalign{&
\int_T^{2T}(S(t+G) - S(t-G))^4\d t \cr&
\ll_\e T^{1+\e}M^{-1}\int_T^{2T}(S^2(t+G) +S^2(t-G))\d t\cr&
\ll_\e T^{1+\e}M^{-1}T^{1/2}\Biggl(\int_T^{2T}\sum_{n>M}d^2(n)n^{-3/2}\cr&
+ \sum_{M<m\ne n\le2M}T^{1/2}d(m)d(n)(mn)^{-3/4}|\sqrt{m}-\sqrt{n}|^{-1}\Biggr)
\cr& \ll_\e T^{3/2+\e}M^{-1}(TM^{-1/2} + T^{1/2}) \ll_\e T^{1+\e}G^2,
\cr}
$$
as requested, since $M \ge TG^{-4/3}$.

\medskip
If
$$
M \le T^{1-\e}G^{-2},\leqno(5.3)
$$
we proceed as follows. First in $S(t\pm G)$ we
replace $(t\pm G)^{1/4}$ by $t^{1/4}$, making a small total error in the process.
Then we note that
$$
\eqalign{&
\cos(4\pi\sqrt{n(t+G)}-\pi/4) - \cos(4\pi\sqrt{n(t-G)}-\pi/4)\cr&
= - 2\sin\left(2\pi\sqrt{n}(\sqrt{t+G} - \sqrt{t-G})\right)
\cos\left(2\pi\sqrt{n}(\sqrt{t+G} + \sqrt{t-G})\right).\cr}
$$
Furthermore, since
$$
\sqrt{t+G} - \sqrt{t-G} = \sqrt{t}\left({G\over t} +
\sum_{j=2}^\infty d_j{\Bigl({G\over t}\Bigr)}^j\right)
\leqno(5.4)
$$
with suitable constants $d_j$, it follows that in view of (5.3) in the series
expansion of
$$
\sin\Bigl(2\pi\sqrt{n}(\sqrt{t+G} - \sqrt{t-G})\Bigr)
$$
the term $2\pi G\sqrt{n/t}$ will dominate in size. Hence if we take sufficiently
many terms in (5.4) the tail of the series will make a negligible contribution,
and we are left with a finite number of integrals, of which the largest
contribution will come from
$$
T\int_T^{2T}{\Biggl|\sum_{M<n\le M'}d(n)n^{-3/4}Gn^{1/2}t^{-1/2}
\exp\Bigl(2\pi i\sqrt{n}(\sqrt{t+G}+\sqrt{t-G})\Bigr)\Biggr|}^4\d t.\leqno(5.5)
$$
Let now $\f(t)\,(\ge 0)$ be a smooth function, supported in $[T/2,\,5T/2]$
and equal to unity in $[T,\,2T]$. Then $\f^{(r)}(t) \ll_r T^{-r}$ for
$r = 0,1,2,...\;$. We have
$$\eqalign{&
\int_T^{2T}|\cdots|^4\d t \le \int_{T/2}^{5T/2}\f(t)|\cdots|^4\d t\cr&
\ll {G^4\over T^2}\int\limits_{T/2}^{5T/2}\f(t)\sum_{k,\ell,m,n\asymp M}
{d(k)d(\ell)d(m)d(n)\over(k\ell mn)^{1/4}}
\exp\Bigl(i\D(\sqrt{t+G} + \sqrt{t-G})\Bigr)\d t,
\cr}
$$
where
$$
\D \;:= \D(k,\ell,m,n) =  \; 2\pi(\sqrt{k}+\sqrt{\ell} - \sqrt{m}-\sqrt{n}\,).
$$
In the last integral we perform a large number of integrations by parts.
During this process the exponential factor will
remain the same, while the integrand will acquire
each time an additional factor of order $\asymp 1/(\D\sqrt{T})$. Hence
the contribution of integer quadruples $(k,\ell,mn)$ for which
$|\D| > T^{\e-1/2}$ will be negligible. The  contribution of the remaining
quadruples is estimated by Lemma 1 (with $k=2, \delta = |\D|T^{\e-1/2}$)
and trivial estimation. In this way it is seen that the expression in (5.5) is
$$
\eqalign{&
\ll_\e T^\e G^4M^{-1}(T^{-1/2}M^{7/2} + M^2)\cr&
= T^{\e-1/2}G^4M^{5/2} + T^\e G^4M\cr& \ll_\e T^{2+\e}G^{-1} + T^{1+\e}G^2
\ll_\e T^{1+\e}G^2
\cr}
$$
for $G\gg T^{1/3}$. It remains to deal with the intermediate range
$$
T^{1-\e}G^{-2} \ll M \ll T^{1+\e}G^{-4/3}.\leqno(5.6)
$$
This is accomplished similarly as in the previous case, by using the
trivial inequality
$$
(S(t+G)-S(t-G))^4 \ll S^4(t+G) + S^4(t-G),
$$
namely by working with two expressions $S(t\pm G)$,
without taking into account the effect of
$t+G$ and $t-G$ combined. We see that the contribution will be, in view of (5.6),
$$\eqalign{&
\ll_\e T^{2+\e}M^{-3}(T^{-1/2}M^{7/2} + M^2)
= T^{3/2+\e}M^{1/2} + T^{2+\e}M^{-1}\cr&
\ll_\e T^{2+\e}G^{-2/3} + T^{1+\e}G^2 \ll_\e T^{1+\e}G^2\cr}
$$
for $T^{3/8} \le G \ll T^{1/2}$, as asserted. This proves
(2.15) and completes the proof of Theorem 4.

\medskip
{\bf Remark 3}. If one had the analogue of (5.1) with $k=2$, namely
the bound $N^\e(N^6\delta + N^3)$ for six square roots, then the above
argument would lead to
$$
\int_T^{2T}(\D(t+G)-\D(t-G))^6\d t \ll_\e T^{1+\e}G^3\quad
(T^{6/13} \le G=G(T) \ll T^{1/2},
$$
which would still be a non-trivial result.
\medskip
In concluding, it may be remarked that one can also obtain another proof of
the important bound
$$
\int_0^T|\zt|^{12}\d t \;\ll_\e\; T^{2+\e}.\leqno(5.7)
$$
This bound is due to D.R. Heath-Brown [4], who had $\log^{17}T$ in place
of $T^\e$, and still represents the sharpest known bound for high
moments of $|\zt|$. Namely in (4.5) we immediately choose $G =
AV^2L^{-2}$ with $t_j = u$ and then integrate, with an additional
smooth weight. Like in the original proof of (5.7) in [4], the sum
$\sum_2(T)$ in Atkinson's formula [1] (or [5, Chapter 15]) for
$E(T)$ will make a negligible contribution, while the range of
summation in $\sum_1(T)$ will be $1 \le n \le T^{1+\e}G^{-2}$. The
technical details are as before, while the function $f(t,n)$ in the
sum $\sum_1(T)$ is neutralized by using a procedure due to M. Jutila
[15, Part II], which was also used in [7]. In this way (5.7) will
eventually follow.


\bigskip\bigskip
\Refs
\bigskip\bigskip

\item{[1]} F.V. Atkinson, The mean value of the Riemann zeta-function,
Acta Math. {\bf81}(1949), 353-376.

\item{[2]} Y. Bugeaud and A. Ivi\'c,
Sums of the error term function in the mean square for $\z(s)$,
to appear, see {\tt arXiv:0707.4275}.

\item{[3]} G. Coppola and S. Salerno, On the symmetry of the divisor
function in almost all short intervals, Acta Arith. {\bf113}(2004),
189-201.

\item{[4]} D.R. Heath-Brown, The twelfth power moment of the Riemann
zeta-function, Quart. J. Math. (Oxford) {\bf29}(1978), 443-462.

\item{[5]} A. Ivi\'c, The Riemann zeta-function, John Wiley \&
Sons, New York, 1985 (2nd ed., Dover, Mineola, N.Y., 2003).

\item{[6]} A. Ivi\'c, The mean values of the Riemann zeta-function,
LNs {\bf 82}, Tata Inst. of Fundamental Research, Bombay (distr. by
Springer Verlag, Berlin etc.), 1991.

\item{[7]}A. Ivi\'c, On the Riemann zeta function and the divisor problem, Central
European Journal of Mathematics 2({\bf4}) (2004),   1-15;
II ibid. 3({\bf2}) (2005), 203-214,
III, subm. to  Ann. Univ. Budapest, Sectio Computatorica, and IV,
Uniform Distribution Theory {\bf1}(2006), 125-135.

\item{[8]} A. Ivi\'c, Some remarks on the moments of $|\zt|$ in short intervals,
to appear in Acta Mathematica Hungarica, see {\tt math.NT/0611427}.

\item{[9]} A. Ivi\'c, On moments of $|\zt|$ in short intervals,
Ramanujan Math. Soc. LNS{\bf2},
The Riemann zeta function and related themes: Papers in honour of
Professor Ramachandra,
 2006, 81-97.

\item{[10]} A. Ivi\'c, On the mean square of the zeta-function
and the divisor problem,
 Annales Acad. Sci. Fennicae Math. {\bf}32(2007), 1-9.

\item{[11]} A. Ivi\'c and P. Sargos, On the higher moments of the
error term in the divisor problem, Illinois J. Math. {\bf81}(2007), 353-377.

\item{[12]} M. Jutila, On the divisor problem for short intervals,
Ann. Univer. Turkuensis Ser. {\bf A}I {\bf186}(1984), 23-30.

\item{[13]} M. Jutila, A method in the theory of exponential sums,
TATA LNs {\bf80}, Springer Verlag, Berlin etc., 1987.

\item{[14]} M. Jutila, Mean value estimates for exponential sums, in
``Number Theory, Ulm 1987", LNM {\bf1380}, Springer Verlag,
Berlin etc., 1989, 120-136.

\item{[15]} M. Jutila, Riemann's zeta-function and the divisor problem,
Arkiv Mat. {\bf21}(1983), 75-96 and II, ibid. {\bf31}(1993), 61-70.

\item{[16]} A. Kaczorowski and A. Perelli, The Selberg class: a
survey, in ``Number Theory in Progress, Proc. Conf. in honour
of A. Schinzel (K. Gy\"ory et al. eds)", de Gruyter,
Berlin, 1999, pp. 953-992.

\item{[17]} M. K\"uhleitner and W.G. Nowak,
The average number of solutions of the
Diophantine quation $u^2+v^2=w^3$ and related arithmetic functions, Acta
Math. Hungarica {\bf104}(2004), 225-240.

\item {[18]} K. Ramachandra and A. Sankaranarayanan, On an asymptotic formula
of Srinivasa Ramanujan, Acta Arith. {\bf109}(2003), 349-357.

\item{[19]}R.~A.~Rankin, Contributions to the theory of Ramanujan's
   function $\tau(n)$ and similar arithmetical functions II. The order
   of the Fourier coefficients of integral modular forms,
    Proc. Cambridge Phil. Soc. {\bf 35}(1939), 357-372.

\item{[20]}R.~A.~Rankin, Modular forms and functions, Cambridge Univ. Press,
Cambridge, 1977.

\item{[21]} O. Robert and P. Sargos, Three-dimensional exponential
sums with monomials, J. reine angew. Math. {\bf591}(2006), 1-20.

\item{[22]}A.~Selberg,  Bemerkungen \"{u}ber eine Dirichletsche Reihe,
    die mit der Theorie der Modulformen nahe verbunden ist,
    Arch. Math. Naturvid. {\bf 43}(1940), 47-50.

\item{[23]} E.C. Titchmarsh, The theory of the Riemann
zeta-function (2nd ed.),  University Press, Oxford, 1986.

\item{[24]} K.-M. Tsang,  Higher power moments of $\D(x)$, $E(t)$
and $P(x)$, Proc. London Math. Soc. (3){\bf 65}(1992), 65-84.

\item{[25]} W. Zhai,  On higher-power moments of $\D(x)$,
Acta Arith.  {\bf112}(2004),
367-395, II ibid. {\bf114}\break (2004), 35-54 and III ibid. {\bf118}(2005), 263-281.

\item{[26]} W. Zhai,  On higher-power moments of $E(t)$,
Acta Arith.  {\bf115}(2004), 329-348.

\vskip2cm
\endRefs

\enddocument

\bye